\documentclass{article}

\pdfoutput=1
\usepackage[T1]{fontenc}
\usepackage{amsfonts}
\usepackage[toc]{appendix}
\usepackage{amsthm}
\usepackage{amsmath}
\usepackage{array}
\usepackage{listings, afterpage, footnote, float}
\usepackage[maxbibnames=99]{biblatex} 

\usepackage{color}
\usepackage{hyperref}
\usepackage{geometry}
\hypersetup{
    pdftitle={A database of constructions of Hadamard matrices},
    colorlinks=true,
    linkcolor=black,
    urlcolor=blue,
    citecolor=black,
    }
\newcommand{\pb}[1]{\sl \color{blue}{#1}}
\newcommand{\rb}[1]{\sf \color{red}{#1}}
\newcommand{\nottt}[1]{\rm \color{black}{#1}}

\theoremstyle{definition}
\newtheorem*{definition}{Definition}
\newtheorem{theorem}{Theorem}

\title{A database of constructions of Hadamard matrices}
\author{Matteo Cati\thanks{Department of Computer Science, Oxford University, UK.\protect\\
\url{matteo.cati@keble.ox.ac.uk}, \url{dmitrii.pasechnik@northwestern.edu}}
\and
Dmitrii V. Pasechnik\footnotemark[1] \footnote{Corresponding author. Current address: Department of Computer Science,
Northwestern University.}}
\date{}

\begin{document}
\maketitle

\begin{abstract}
Hadamard matrices of order $n$ are conjectured to exist whenever $n$ is $1$,
$2$, or a multiple of $4$; a similar conjecture exists for skew
Hadamard matrices.  We provide constructions covering orders $\le 1208$
of all known Hadamard and skew Hadamard matrices in the open-source
software SageMath. This allowed us to verify the correctness of results
given in the literature.  Within this range, just one order, $292$, of
a skew Hadamard matrix claimed to have a known construction, required a fix.

We also produce the up to date tables, for $n \le 2999$ (resp. $n\le 999$ for
skew case), of the minimum exponents $m$ such that a (skew) Hadamard
matrix of order $2^m n$ is known, improving over 100 entries in the
previously published sources. We explain how tables' entries are
related to Riesel numbers. As a by-product of the latter, we show that
the Paley constructions of (skew-)Hadamard matrices do not work for the
order $2^m 509203$, for any $m$.
\end{abstract}

\section{Introduction}
\label{intro-chapter}
A $n \times n$ matrix is called a \textit{Hadamard matrix} of order $n$ if its entries are all $\pm 1$, and the rows are mutually orthogonal. Equivalently, they are $\pm 1$ matrices which satisfy the equation:
\begin{displaymath}
    HH^\top = H^\top H = nI,
\end{displaymath}
where $I$ is the identity matrix.
Matrices of this type were described for the first time by J.~Sylvester \cite{Sylvester1867}, and were studied further by J.~Hadamard \cite{Hadamard1893}.

Additionally, a Hadamard matrix is  \textit{skew} if it has the form $H=S+I$, where $S$ is skew-symmetric: $S^\top = -S$.

These matrices are used in many practical applications, such as data compression \cite{tropp2011},
image analysis, signal processing \cite{RYJH1997}, statistics \cite{Rsystem, Rhadamard},
quantum computing \cite{Suksmono2022}.
However, given that many constructions exist, each covering only some orders, finding the correct one for a specific order is often difficult.

The purpose of this paper is checking, detailing, and ``digitally preserving'' a number of constructions, and tables specifying known constructions and orders. For the latter in particular, we updated about 100 entries in two tables in J.~Seberry's and M.~Yamada's monograph \cite{Seberry2020}.
We provide all the necessary constructions in the open-source mathematical software SageMath \cite{sagemath}. While it is easy to create and manipulate Hadamard matrices in SageMath
installed on one's computer, the corresponding data can also be created and extracted online (without a local installation of SageMath), using a web interface, as explained in  Sect.~\ref{runsage}.

The relevant parts of SageMath evolved over the past 15 years, and comprise roughly 15,000 lines of well-documented Python code+tests, 
with a number of constructions of Hadamard matrices and related combinatorial constructions, such
as difference sets, finite geometries, etc.
already present in SageMath in 2014. Before the coding push in 2023 which led to this paper,
SageMath had a number of constructions of regular symmetric Hadamard
matrices with constant diagonal
(\href{https://doc.sagemath.org/html/en/reference/combinat/sage/combinat/matrices/hadamard_matrix.html#sage.combinat.matrices.hadamard_matrix.regular_symmetric_hadamard_matrix_with_constant_diagonal}{RSHCDs}).
An RSHCD corresponds to a strongly regular graph (SRG); they were implemented as part of
N.~Cohen's and D.~Pasechnik's work \cite{Cohen2017}, which provides
constructions for SRGs with at most 1300 vertices and sets of parameters for which existence is known,
as in A.E.~Brouwer \cite{aebdb}.

In 2023 we added a number of constructions of Hadamard and skew Hadamard matrices, so that
all known orders $\le 1208$ are now available. A similar programme for symmetric
Hadamard matrices (cf. e.g. \cite{ddj2022}) is planned. The first order for which a
Hadamard matrix exists, but is not implemented in Sagemath is 1212 (the order 1212 example is
due to D.~Djoković and I.S.~Kotsireas \cite{Djokovic2016}),  while the first missing skew Hadamard matrix
is the one of order 1252 (due to J.~Seberry \cite{Seberry2017}, Table 4.5).

We note that currently, cf. \cite{Djokovic2016} or Table~\ref{hadamard-known-list-table} in Appendix,
the only orders in this range for which Hadamard matrices are unknown are 668, 716, 892 and 1132.
For the skew case, 41 are unknown, cf. Table~\ref{skew-hadamard-known-list-table}.\footnote{There are first 36
entries, from 89 to 299, of which 5 account for two matrices, e.g. 89 accounts for orders 356 and 712.}
%

In addition to making known matrices available in SageMath, implementing known
constructions allowed us to verify their correctness. In particular, we
were unable to locate in the literature a (correct) construction for a skew Hadamard matrix of order 292.
In Sect.~\ref{missing292} we describe a new construction for it, which was kindly
provided to us by D.~Djoković \cite{Djokovic2023pers}.

Constructions use a variety of algebraic combinatorics objects, 
briefly discussed in Sect.~\ref{background-chapter}.
Not all of these were
present in SageMath---in these cases, we implemented the missing constructions.

Sect.~\ref{design-chapter} briefly explains the
SageMath interface to the constructions, and provides
hints on running SageMath.
Sect.'s~\ref{hadamard-chapter} and \ref{skewmats}
describe the implemented constructions of general
and skew-Hadamard matrices.
In Sect.~\ref{knownhm} we describe tables, for odd
$n \le 2999$ (respectively $n\le 999$ for skew matrices),
of the minimum exponents $m$ such that a (skew) Hadamard matrix of order $2^m n$ is known,
filling gaps in \cite{Seberry2020}, Tables A.17 and 9.2. The corrections are based on our
implementation of matrix constructions, as well as on the observation
that in the tables' range, there should be no entries for which $m$ is unknown,
thanks to results on \emph{Riesel numbers}, cf. Sect.~\ref{sect:riesel}.
One remarkable by-product of the latter
is the following result: there is no Hadamard matrix of a Paley type for order $2^m\times 509203$, for any $m$,
cf. Theorem~\ref{thm:riesel}.
Finally, Appendix~\ref{appendix} contains descriptions of the tables
of known matrices, and the tables themselves.

\subsection*{Related work}

Implementing constructions is a first step towards a formal verification of them, using a proof assistant. Work towards
such a verification has been done by Lu-Ming Zhang (supervised by the 2nd author),
for constructions of Hadamard matrices covering all the possible orders up to 112 \cite{luming2021}. The latter
demonstrated challenges in formal checking of constructions of Hadamard matrices. Currently, proof assistants are
severely deficient in doing even slightly nontrivial calculations, and closing this lacuna is an active research area.  Previous efforts to combine proof assistants and computer algebra systems (CASs, for short),
e.g.   
\cite{bauer98,davenport2023,lewis2021}
have demonstrated that one either has to trust the CAS and thereby significantly lower correctness guarantees, or severely hamper the power of CAS by insisting that it produce verifiable certificates (which CASs are not designed to do at present).

S.~London in his PhD thesis \cite{MR3187587} outlined constructions of matrices
for all known orders not exceeding $1000$, although he did not provide an
implementation. We took from [loc.cit.] data for a number of small
Williamson type Hadamard matrices. The remainder of our work is independent of [loc.cit.]---in particular as
it did not touch the skew case.

Riesel numbers, in conjunction with the asymptotic existence of
Hadamard matrices, were used by W.~de~Launey \cite{DELAUNEY20091002}. Our application is
different, see Sect~\ref{sect:riesel} for more details.

\subsubsection*{Other online databases and computer systems}
N.J.A. Sloane maintains an online database of Hadamard matrices of order $\leq 256$ \cite{sloantables}.
A bigger database \cite{hadgen} by V.K.~Gupta, R.~Parsad and  A.~Dhandapani appears to be permanently
offline.\footnote{The website is still there, but attempts to query a matrix return ``Internal server error''.} 

Magma \cite{MR1484478}, a commercial computer algebra system,
provides extensive databases of Hadamard matrices, general and skew, up to order 256 \cite{magmatables}.
Magma's documentation does not list functions to generate Hadamard matrices of arbitrary size bigger than 256.

Statistics computing open-source system R \cite{Rsystem}
provides a function {\tt hadamard(n)} which, given $n$, constructs a Hadamard matrix of order greater than
$n$ \cite{Rhadamard}. It only uses Paley (Sect.~\ref{paleyconst})
and Sylvester (Sect.~\ref{doublingconst})
constructions. E.g. it is unable to construct
a matrix of order 92---{\tt hadamard(90)} will output a matrix of order 96.

As we see, our project covers considerably more ground than these listed above. Indeed, only \cite{Rhadamard}
provides implementations of Hadamard matrix constructions, but it only uses two most popular constructions.
And the databases only go up to order 256, compared to our 1200.

\section{Additional Definitions}
\label{background-chapter}
Several mathematical objects were needed in order to implement the necessary constructions of Hadamard matrices. Here we give a brief overview of such objects:
T-sequences and matrices; difference sets;
Williamson, and good matrices; Paley and Paley-like
(di)graphs.

\subsection{T-sequences and T-matrices}\label{sect:Tmat}
Let $N_X(j)$ denote the \emph{non-periodic autocorrelation} of an $n$-sequence 
$X=(x_1,\dots,x_n)$, as in \eqref{NPAF}, which is (1.14) in \cite{Seberry2020}.
\begin{equation}
    \label{NPAF}
    N_X(j) =\sum_{i=1}^{n-j}x_i x_{i+j},\qquad 0\leq j\leq n-1.
\end{equation}
T-sequences are defined by \cite{Kharaghani2005} as:
\begin{definition}
    Four $(0, \pm 1)$ sequences $A, B, C, D$ of length $n$ are called T-sequences if
\begin{displaymath}
    N_A(j) + N_B(j) + N_C(j) + N_D(j) = 0  \text{\: for \:}  j\ge 1
\end{displaymath}
    and in each position exactly one of the entries of $A, B, C, D$ is nonzero.
\end{definition}

To implement T-sequences in SageMath, we used the construction detailed in
H.~Kha\-ra\-ghani \cite{Kharaghani2005}, as well as Theorem 7.7 of
\cite{Seberry2017}.\footnote{\url{https://doc.sagemath.org/html/en/reference/combinat/sage/combinat/t_sequences.html}
gives implementation details and examples.}

These sequences are needed to create T-matrices, which are defined in \cite{Seberry2017} as follows:

\begin{definition}
    Four circulant $(0, \pm 1)$ matrices $X_i$, $i = 1,2,3,4$, of order $n$ which are non-zero for each of the $n^2$ entries for exactly one $i$, and which satisfy
\begin{displaymath}
    \sum_{i=1}^4X_iX_i^\top = nI
\end{displaymath}
will be called T-matrices of order $n$.
\end{definition}

If we use T-sequences of length $n$ as the first rows of four circulant matrices $X_1$, $X_2$, $X_3$, $X_4$,
we obtain four T-matrices. Additionally, some T-matrices have been computed directly in \cite{Sawade1985}.
T-matrices are used in Theorem~\ref{CooperWallisConstruction-section}.

\subsection{Supplementary and Complementary Difference Sets}\label{sect:SDS}
The following is Definition 4.3 from \cite{Seberry2017}.
\begin{definition}
    Let $S_1, S_2, ..., S_n$ be subsets of $G$, an additive abelian group
of order $v$. Let $|S_i| = k_i$. If the equation $g = r - s$, $r, s \in S_i$ has exactly $\lambda$ solutions for each non-zero element $g$ of $G$, then $S_1, S_2, ..., S_n$ will be called $n - \{v;k_1, k_2, ...,k_n;\lambda\}$ supplementary difference sets (or SDS). If $k_1 = k_2 = ... = k_n = k$, we call it a $n - \{v;k;\lambda\}$ SDS.
\end{definition}

We have implemented a number of constructions of SDSs.\footnote{See \url{https://doc.sagemath.org/html/en/reference/combinat/sage/combinat/designs/difference_family.html\#sage.combinat.designs.difference_family.supplementary_difference_set}
for details and examples.}
Specifically, we use results from small order SDSs from \cite{Djokovic1992a, Djokovic1992b, Djokovic1992c, Djokovic1994a,
Djokovic2008a, Djokovic2008b, DGK2014, Djokovic2016}.

We also need SDSs with the additional properties
\begin{enumerate}
    \item $S_1 \cap -S_1 = \emptyset$;
    \item $S_1 \cup -S_1 = G\setminus \{0\}$.
\end{enumerate}
In this case, we say that $S_1$ is \emph{skew}, and the sets are called \emph{skew} SDSs.
\footnote{\url{https://doc.sagemath.org/html/en/reference/combinat/sage/combinat/designs/difference_family.html\#sage.combinat.designs.difference_family.skew_supplementary_difference_set}
gives implementation details and examples.}
SDSs with parameters $2-\{2m+1; m, m; m-1\}$
are called \emph{complementary difference sets} of size $m$, cf. G.~Szekeres \cite{Szekeres1969}.
We implement constructions of such sets from  \cite{Szekeres1969, Szekeres1971}.
SDSs are used in constructions from Sect.~\ref{gsarray}.

\subsection{Williamson and Good Matrices}\label{sect:good}
These are particular kinds of matrices used to generate Hadamard matrices
in 

\begin{definition}
    Let $A, B, C, D$ be $n\times n$ symmetric circulant matrices,
satisfying the condition:
    \begin{displaymath}
        A^2+B^2+C^2+D^2 = 4nI
    \end{displaymath}
    Then, these matrices are called Williamson matrices.
\end{definition}

Williamson matrices are available for infinitely many values of $n$, cf. e.g. M.~Hall \cite{Hall1988}---
this construction is not yet implemented in SageMath. Data for small cases, up to size 63, is available.\footnote{See 
\url{https://doc.sagemath.org/html/en/reference/combinat/sage/combinat/matrices/hadamard_matrix.html\#sage.combinat.matrices.hadamard_matrix.williamson_type_quadruples_smallcases} for details.}
They are used in the construction in Sect.~\ref{sect:willi}.

\begin{definition}
    Four $\pm 1$ matrices $A, B, C, D$ of order $n$ (odd) with the properties:
\begin{itemize}
    \item $MN^\top = NM^\top$ for $M, N \in \{A, B, C, D\}$;
    \item $(A-I)^\top = -(A-I)$, $B^\top = B$, $C^\top = C$, $D^\top = D$; 
    \item $AA^\top + BB^\top + CC^\top + DD^\top = 4nI$.
\end{itemize}
are called good matrices.
\end{definition}

Examples of good matrices of order $n=1, 3, ..., 31$ are listed in \cite{Szekeres1987}.\footnote{See
\url{https://doc.sagemath.org/html/en/reference/combinat/sage/combinat/matrices/hadamard_matrix.html\#sage.combinat.matrices.hadamard_matrix.skew_hadamard_matrix_from_good_matrices} for details.}
They are used in constructions in Sect.~\ref{skewmats}. 

\subsection{Paley (di)graphs and generalisations}\label{sect:paleydigraphs}
Let $q$ be an odd prime power. Define a $q\times q$ matrix $\Pi=(P_{ij})$, with rows and columns indexed by the elements
of the finite field $\mathbb{F}_q$, as follows.
\begin{equation*}
    \Pi_{ij}=\begin{cases}
        \phantom{-}1 & \text{if $i-j$ a square in $\mathbb{F}_q$}\\
        -1 & \text{otherwise}
    \end{cases}
\end{equation*}
The matrix $\Pi$ can be viewed as the adjacency  matrix of a
graph on $\mathbb{F}_q$, which is undirected if $q\equiv 1\pmod 4$,\footnote{See
\url{https://doc.sagemath.org/html/en/reference/graphs/sage/graphs/graph_generators.html\#sage.graphs.graph_generators.GraphGenerators.PaleyGraph} for details.}
and a tournament (a directed graph with arcs present between 
any unordered pair of vertices) otherwise.\footnote{See \url{https://doc.sagemath.org/html/en/reference/graphs/sage/graphs/digraph_generators.html\#sage.graphs.digraph_generators.DiGraphGenerators.Paley} for details.}
Here 1 indicates adjacency, and loops are ignored.

If $q\equiv 3\pmod 4$ then $\Pi$ is skew-symmetric, and
$\Pi\Pi^\top=(q+1)I-J$,
where $J$ denotes the all-1 matrix. 
If  $q\equiv 1\pmod 4$ then $\Pi$ is symmetric, and
$\Pi^2=2\Pi+(q-1)I-J$.

More generally, there are constructions of (di)graphs with these algebraic properties, but $q$ not necessarily a prime power,
e.g. R.~Mathon \cite{Mat1978} constructs such graphs of order $pq^2$, with $p,q$ odd prime powers.
Matrices just described are used in the constructions in Sect.~\ref{paleyconst}.

\section{Program Design}
\label{design-chapter}

In SageMath, all Hadamard matrices are available via a single function,
\href{https://doc.sagemath.org/html/en/reference/combinat/sage/combinat/matrices/hadamard_matrix.html}{\tt hadamard\_matrix} from the module {\tt sage.combinat.matrices.hadamard\_matrix}.

Similarly, all skew Hadamard matrices can be obtained using
\href{https://doc.sagemath.org/html/en/reference/combinat/sage/combinat/matrices/hadamard_matrix.html#sage.combinat.matrices.hadamard_matrix.skew_hadamard_matrix}{\tt skew\_hadamard\_matrix} from the same module.

\begin{lstlisting}[language=python, basicstyle={\ttfamily}]
    sage: from sage.combinat.matrices.hadamard_matrix
    ....: import hadamard_matrix, skew_hadamard_matrix
    sage: H = hadamard_matrix(4)
    sage: H
    [ 1  1  1  1]
    [ 1 -1  1 -1]
    [ 1  1 -1 -1]
    [ 1 -1 -1  1]
    sage: H2 = skew_hadamard_matrix(20)
    sage: H2
    20 x 20 dense matrix over Integer Ring
\end{lstlisting}
The first parameter of these functions is the order of the Hadamard (or skew Hadamard) matrix that should be returned. Furthermore, they have optional boolean parameters:
\begin{itemize}
    \item {\tt existence} is {\tt False} by default. When it is set to {\tt True}, the function will return either {\tt True}, {\tt False} or {\tt Unknown} instead of a matrix. In particular, it will return {\tt True} if the matrix of the given order can be constructed, {\tt Unknown} if a construction for that matrix is not yet present in SageMath (or it is not known), and {\tt False} if the matrix does not exist (i.e. $n$ is not 1, 2 or a multiple of 4).
    \item {\tt check} is {\tt True} by default. Normally, the functions will perform a check on the generated matrix before returning it, to ensure that it satisfies the requirements. This check will be skipped if the parameter is set to {\tt False}. Note that the value will have no effect if {\tt existence} is {\tt True}.
    \item {\tt construction\_name} is {\tt False} by default, and has no effect is {\tt existence} is  {\tt False}.
    If both {\tt construction\_name} and {\tt existence} are {\tt True}, a short info about the construction will be reported.
\end{itemize}
\begin{lstlisting}[language=python, basicstyle={\ttfamily}]
    sage: hadamard_matrix(12, existence=True, construction_name=True)
    'paleyII 12'
    sage: hadamard_matrix(668, existence=True)
    Unknown
    sage: hadamard_matrix(9, existence=True)
    False   
\end{lstlisting}
The function {\tt skew\_hadamard\_matrix} has one more optional parameter, {\tt skew\_normalize}. When this is set to {\tt True} (the default value) the skew Hadamard matrix returned by the function will be normalized, i.e. the first row will contain all $+1$.

Furthermore, each construction added to SageMath is implemented as a separate function, all contained in the module {\tt sage.combinat.matrices.hadamard\_matrix}. A complete description of these functions can be found in the SageMath Manual \cite{sagemath}.

\begin{lstlisting}[language=python, basicstyle={\ttfamily}]
    sage: from sage.combinat.matrices.hadamard_matrix 
    ....: import hadamard_matrix_paleyI
    sage: H = hadamard_matrix_paleyI(28)
    sage: H
    28 x 28 dense matrix over Integer Ring
\end{lstlisting}

\subsection{Running SageMath, offline and online}\label{runsage}
SageMath is a Python-based system, freely available and relatively
easy to \href{https://doc.sagemath.org/html/en/installation/index.html}{install}. However, one can also run SageMath in a web
browser, using \href{https://cocalc.com}{CoCalc} (requires creating a free account),
or \href{https://sagecell.sagemath.org/}{SageCell} (just open the link
and execute commands in a browser window), meaning that
in effect the result of our work can be viewed as an online
database.

\section{Constructions for Hadamard Matrices}
\label{hadamard-chapter}
Here we give a brief overview of the constructions that we used. For a more detailed description of how these were implemented, see \cite{cati2023}.

\subsection{Paley's Constructions}\label{paleyconst}
In 1933, Paley \cite{Paley1933} discovered two constructions of Hadamard matrices. 
They use the matrix $\Pi$ from Sect.~\ref{sect:paleydigraphs}
and are as follows, with $\overline{1}$ denoting the column vector
of all 1s.

\begin{theorem}
\label{paleyI-theo}
Let $q \equiv 3 \pmod 4$, and a (generalised) Paley digraph, cf. Sect.~\ref{sect:paleydigraphs},
of order $q$ exists (e.g. $q$ is a prime power).
Then there is a skew Hadamard matrix $H$ of order $q+1$,
obtained from $\Pi$ by taking
$H=\begin{pmatrix}
    1 & \overline{1}^\top\\
    -\overline{1} & \Pi
\end{pmatrix}$. \qed
\end{theorem}

\begin{theorem}\label{paleyII-theo}
Let $m = 2(q+1)$, for $q \equiv 1 \pmod 4$
 a (generalised) Paley graph, cf. Sect.~\ref{sect:paleydigraphs},
of order $q$ exists (e.g. $q$ is a prime power). 
Then there is a symmetric Hadamard matrix $H$ of order $m$, obtained
from $\Pi$ by taking \\
$H=I_{q+1}\otimes
\begin{pmatrix}
1&-1\\-1&-1    
\end{pmatrix}+\begin{pmatrix}
    0 & \overline{1}^\top\\
    \overline{1} & \Pi-I_q
\end{pmatrix}\otimes
\begin{pmatrix}
    1&1\\1&-1    
    \end{pmatrix}$.
\qed
\end{theorem}

\subsection{Doubling Constructions}\label{doublingconst}
Nearly 60 years before Paley, Sylvester \cite{Sylvester1867} had already proved that if $n$ is the order of a Hadamard matrix, then there exists a Hadamard matrix of order $2^tn$ for all values of $t\ge0$.

\begin{theorem}
Let $H$ be a Hadamard matrix of order $n$. Then the matrix $H' = \begin{pmatrix}
H & H\\
H & -H
\end{pmatrix}$ is a Hadamard matrix of order $2n$. \qed
\end{theorem}

More recently, Seberry \cite{Wallis1971} described a similar construction for skew Hadamard matrices:
\begin{theorem}
    Suppose $H_n = S + I_n$ is a skew Hadamard matrix of order $n$. Then

    \[
    H_{2n} = 
    \begin{pmatrix}
    S+I_n & S+I_n\\
    S-I_n & -S+I_n
    \end{pmatrix}
    \]

    is a skew Hadamard matrix of order $2n$. \qed
\end{theorem}

\subsection{Williamson Construction} \label{sect:willi}
Another construction for Hadamard matrices is due to Williamson, see \cite{Hall1988}, Sect. 14.2. He proved that the matrix $H$ given by 

\begin{displaymath}
H = \begin{pmatrix}
A & -B & -C & -D \\
B & A & D & -C \\
C & -D & A & B \\
D & C & -B & A
\end{pmatrix}
\end{displaymath}

is a Hadamard matrix of order $4n$, if $A, B, C, D$ are Williamson matrices.

\subsection{Goethals-Seidel Array}\label{gsarray}
Goethals and Seidel \cite{Goethals1970} discovered that, given four $\pm 1$ matrices $A, B, C, D$ of order $n$ such that
\begin{displaymath}
    AA^\top+BB^\top+CC^\top+DD^\top = 4nI,
\end{displaymath}
a Hadamard matrix of order $4n$ can be constructed by plugging them into the Goethals-Seidel array
\begin{displaymath}
    GS(A, B, C, D) = \begin{pmatrix}
        A & BR & CR & DR \\
        -BR & A & D^\top R & -C^\top R \\
        -CR & -D^\top R & A & B^\top R \\
        -DR & C^\top R & -B^\top R &  A
    \end{pmatrix}
\end{displaymath}
Here, $R$ is the $n\times n$ permutation matrix with all ones on the anti-diagonal.

Furthermore, if $A = S+I$ with $S$ a skew-symmetric matrix, the Goethals-Seidel array will give a skew Hadamard matrix. 

The Goethals-Seidel array has been used extensively to construct Hadamard matrices of order $4v$ from SDSs with parameters $4-\{n; k_1, k_2, k_3, k_4; n-\sum_{i=1}^4k_i\}$ \cite{Djokovic2008}.
Additionally, a skew SDS can be used to obtain a skew Hadamard matrix.

J.~Cooper and J.~Wallis \cite{Cooper1972} showed a construction which uses T-matrices and a generalisation of Williamson matrices:
\begin{theorem}\label{CooperWallisConstruction-section}
Suppose there exist four  T-matrices of order $n$. Furthermore, suppose there exist four $\pm 1$ matrices $A$, $B$, $C$, $D$ of order $w$ which pairwise satisfy $MN^\top$ = $NM^\top$ and for which
\begin{displaymath}
AA^\top + BB^\top + CC^ \top + DD^ \top = 4wI
\end{displaymath}
Then there exists a Hadamard matrix of order $4nw$ .\qed
\end{theorem}

In \cite{Miyamoto1991}, Miyamoto provided a construction for an infinite series of Hadamard matrices:
\begin{theorem}\label{miyamoto}
    Let $q$ be a prime power and $q \equiv 1 \mod 4$. If there is a Hadamard matrix of order $q-1$, then there is a Hadamard matrix of order $4q$. \qed
\end{theorem}

Spence \cite{Spence1975} proposed a construction that uses supplementary difference sets with parameters $4-\{2v; v, v, v, v+1; 2v\}$:
\begin{theorem}\label{sect:spence}
    If there exist $4-\{2v; v, v, v, v+1; 2v\}$ supplementary difference sets in the cyclic group of order $2v$ then there exists a Hadamard matrix of order $4(2v+1)$. \qed
\end{theorem}

\section{Constructions for skew Hadamard Matrices}\label{skewmats}
Finally, we give a short description of a few constructions that can be used to obtain more skew Hadamard matrices.

A first construction is detailed by Szekeres \cite{Szekeres1987}.
\begin{theorem}
    Given four Good matrices $A, B, C, D$ of order $n$, the following is a skew Hadamard matrix of order $4n$:

\begin{displaymath}
H = \begin{pmatrix}
A & B & C & D \\
-B & A & D & -C \\
-C & -D & A & B \\
-D & C & -B & A
\end{pmatrix}
\end{displaymath}
\end{theorem}

A construction due to Williamson is described in Lemma 14.1.5 of \cite{Hall1988}.
\begin{theorem}\label{lemma1415hall}
If there is a skew Hadamard matrix of order $n$, and for a prime $p$ and $k\geq 1$ one has $p^k+1\equiv 0 \pmod 4$, 
then there is a skew Hadamard matrix of order $n(p^k+1)$. \qed
\end{theorem}

Furthermore, Blatt and Szekeres showed in \cite{Blatt1969} that
\begin{theorem}
 If for some abelian group G of order $2m + 1$ there exists a pair of complementary difference sets $A, B$, then there exists a skew Hadamard matrix of order $4(m + 1)$. \qed
\end{theorem}

A different construction which also uses complementary difference sets was given by Spence \cite{Spence1975b}:

\begin{theorem}\label{thm:spence75b}
    If there exists a cyclic projective plane of order $q$ and two complementary difference sets in a cyclic group of order $1 + q + q^2$, then there exists a skew Hadamard matrix of the Goethals-Seidel type of order $4(1 + q + q^2)$. \qed
\end{theorem}

Lastly, skew Hadamard matrices can be obtained from orthogonal designs \cite{Seberry1978}:

\begin{theorem}\label{thm:amicable}
    Suppose there is an orthogonal design of type $(1,m,mn-m-1)$ in order $mn$. Suppose $n$ is the order of amicable designs of types $((1, n-1); (n))$. Then there is a skew Hadamard matrix of order $mn(n-1)$. \qed
\end{theorem}

\subsection{Skew Hadamard Matrix of Order 292}\label{missing292}
According to multiple papers \cite{Koukouvinos2008, Seberry1978, Seberry1991}, a skew Hadamard matrix of order 292 was known since 1978. However, the only paper which we were able to find describing a construction for it was a paper by Djoković \cite{Djokovic2010}, which incorrectly cites the construction for a non-skew Hadamard matrix.

We contacted the author of \cite{Djokovic2010}, who acknowledged that the paper did not contain the correct reference, and helped us by constructing himself the skew SDS that can be used to obtain the skew Hadamard matrix
\cite{Djokovic2023pers}.

For completeness, below is the (implemented)
construction of such supplementary difference sets.

Let $G$ be the multiplicative group of the nonzero residue classes modulo $n=73$, and let $H = \{1,2,4,8,16,32, 37, 55, 64\}$ be a subgroup of order 9.

Then, we define the 8 cosets $\alpha_0, \alpha_1, ..., \alpha_7$ of $H$ in $G$ such that $ \alpha_{2i+1} = -1 \cdot\alpha_{2i}$ and:
\begin{align*}
    \alpha_0 &= H &  \alpha_2 &= 3H & \alpha_4 &= 5H & \alpha_6 &= 11H
\end{align*}

Now, consider the four index sets:
\begin{align*}
    J_1 &= \{1,3,4,6\} &
    J_2 &= \{3,5,6,7\} &
    J_3 &= \{1,4,5,7\} &
    J_4 &= \{0, 2, 7\}
\end{align*}
The skew SDS are given by:
\begin{align*}
    S_1 &= \bigcup_{j \in J_1} \alpha_j &
    S_2 &= \bigcup_{j \in J_2} \alpha_j &
    S_3 &= \bigcup_{j \in J_3} \alpha_j &
    S_4 &= \{0\} \cup \bigcup_{j \in J_4} \alpha_j
\end{align*}

These sets are used to create a skew Hadamard matrix using the construction from Sect.~\ref{gsarray}. 
One obtains four circulant matrices $A_1, A_2, A_3, A_4$ from them \cite{Djokovic2023}. 
Consider e.g. the matrix $A_1$. 
The first row will be $(a_0, a_1, ..., a_{72})$ where $a_i = -1$ if $i \in S_1$ and  $a_i = 1$ otherwise. 
When plugged into the Goethals-Seidel array, these matrices generate a skew Hadamard matrix of order $292$.

\section{Known Hadamard Matrices}\label{knownhm}

J.~Seberry and M.~Yamada published (cf. \cite{Seberry2020}, Table A.17) a table
which contains for every $n \le 2999$  the minimum exponent $m$ such that a
Hadamard matrix of order $2^m n$ is known, as well as a similar table (cf.
\cite{Seberry2020}, Table 9.2 pp.198-200) for skew Hadamard matrices (with $n
\le 999$).
Tables \ref{hadamard-known-list-table} and \ref{skew-hadamard-known-list-table} in Appendix
contain an updated version of these lists.
There for every odd $n$ we write $n(m)$ to indicate that $m$ is the smallest
exponent such that a Hadamard matrix of order $2^m n$ is known. If $n$ is not
present, it means that $m = 2$.

Table \ref{hadamard-known-list-table} improves 22 entries with $m>2$ in Table A.17
\cite{Seberry2020} - for these we are able to construct examples with $m=2$
using methods implemented - specifically, ones in Sect.~\ref{gsarray}, Theorem~\ref{CooperWallisConstruction-section},
Theorem~\ref{miyamoto}, and Sect.~\ref{skewmats}. Furthermore, we note that Table A.17 \cite{Seberry2020} contains a typo, as the smallest exponent
for $179$ should be $m=3$, not $m=2$. We corrected this in Table \ref{hadamard-known-list-table}.

The situation with Table \ref{skew-hadamard-known-list-table} is more interesting.
In particular, $m$'s listed as unknown (indicated by ``$.$'') in Table 9.2 \cite{Seberry2020} are now filled in.
In more detail, we improve 25 entries with $m>2$ in Table 9.2 \cite{Seberry2020} to $m=2$,
using methods implemented. Additionally, a skew Hadamard matrix of order $4\cdot 397$ has been recently found by Djoković \cite{djokovic2024two}.
For further 42 entries, the values of $m$ are improved downward, although not to $m=2$, again using methods implemented.
Lastly, 14 entries (191, 449, 491, 577, 643,  659, 709, 733, 739, 829, 877, 929, 933, 954) with $m$
listed as unknown in \cite{Seberry2020}, are filled in
via Paley construction and a bit of theory, see Sect.~\ref{sect:riesel},
using the table from \cite{oeisA040081}. Note that a rather big value of $m$ for $n=659$ is hard to even verify
without special tools for primality testing.

Sect.~7 of \cite{Miyamoto1991} claims existence of Hadamard matrices of a number of orders listed (i.e., unknown) in
Tables \ref{hadamard-known-list-table}
and \ref{skew-hadamard-known-list-table}. E.g. order
$4\times 515$ is claimed to be known. However, due to missing details, we were
unable to verify these constructions.

\subsection{Riesel numbers}\label{sect:riesel}
While working on updating Tables~\ref{hadamard-known-list-table} and \ref{skew-hadamard-known-list-table},
we observed that the unfilled prime power entries $p^k$, $k\geq 1$, within the tables' ranges, may be
filled in by using an appropriate Paley
construction. This gives a natural question: how
powerful are the latter, will they work regardless
of the value of $p^k$ ? 
It turns out that the answer is no.

An odd natural number $k$ is called a \emph{Riesel number} if $2^m k -1$ is composite for all $m\geq 0$,
see \cite{Riesel1956, RGuy}. Specifically,
in \cite{Riesel1956} H.~Riesel constructed an infinite family of such numbers,
starting from $509203$, which is believed to be the smallest such number. Computational work
has been carried out to verify the latter, and the smallest number with unknown status is $23699$ (May 2024),
see \cite{Rieselwiki}.
Riesel numbers, in connection to asymptotic existence of Hadamard matrices, have been used by W.~de~Launey \cite{DELAUNEY20091002},
and the first part of Theorem~\ref{thm:riesel} is implicit there. The second part appears to be new.
\begin{theorem}\label{thm:riesel}
Let $k<23669$ be odd, or more generally, let $k$ be any non-Riesel odd number.
Then for some $m>1$ there exists a skew Hadamard matrix of Paley type,
of order $2^m k$. As well, there are no Hadamard matrices of Paley types
(neither symmetric nor skew-symmetric) for order $2^m 509203$, for any $m$.
\end{theorem}
\begin{proof}
The only part which is not immediate from the results mentioned above is why
there is no prime $p$ such that $2^m 509203-1=p^t$, for some $t>0$. Note that
$t=1$ is ruled impossible by \cite{Riesel1956}. To sort out the case $t>1$ we need
more technique used in [loc.cit.], as communicated to us by Max Alekseyev \cite{MORiesel};
this argument applies to the case $t=1$, as well.

Denote $k=509203$. Any number of the form $2^m k-1=p^t$ is divisible by $p\in\{3,5,7,13,17,241\}$,
as shown in \cite{Riesel1956}.
Thus our $p$ must come from this set. By observing that $2^m k-1 \equiv p^t\pmod{M}$
has no solution for an appropriate $M\in\mathbb{N}$, we obtain our result.
Namely, the values of $M$ to achieve this are as follows:
 $85$ for $p=3$;
 $31$ for $p=5$; 
 $57$ for $p=7$;
 $105$ for $p=13$;
 $15$ for $p=17$ and $p=241$.
 Thus it is not possible to construct a skew-Hadamard matrix of order $2^m 509203$ using the construction
 from Theorem~\ref{paleyI-theo}.

 It is also not possible to construct a symmetric Hadamard matrix of this order using the construction
 from Theorem~\ref{paleyII-theo}. Indeed,
 for this we need $\frac{2^m 509203-2}{2}=2^{m-1} 509203-1$ to be a prime power, which we have shown to be impossible.
\end{proof}

In general, different construction methods are needed
to get matrices of order $2^\tau
p^k$. Research on such methods is surveyed 
in Chapter~11 of \cite{Seberry2020}---
it is shown that for any odd $q$ there exists
$\tau_0$ so that a Hadamard matrix of order $2^\tau q$ exists for any $\tau\geq\tau_0$, and results on the best
$\tau_0$ are surveyed.

\subsection*{Acknowledgements}
The authors thank Prof.~Dragomir~\v{Z}.~Djokovi\'{c} for help with
constructing a skew Hadamard matrix of order 292, and many
helpful remarks.
The authors thank Prof.~Max~Alexeyev for \cite{MORiesel},
crucial in Sect.~\ref{sect:riesel}.
The authors thank the anonymous referees for many
helpful comments and suggestions.

\printbibliography[heading=bibintoc]

@book{Hall1988,
  author    = {Marshall Hall},
  address   = {Hoboken, NJ, USA},
  doi       = {10.1002/9781118032862},
  isbn      = {978-1118032862},
  publisher = {John Wiley \& Sons, Inc.},
  title     = {Combinatorial Theory},
  year      = {1988}
}

@book{RYJH1997,
  title     = {Hadamard Matrix Analysis and Synthesis,
               With Applications to Communications and Signal/Image Processing},
  author    = {R. K. Rao Yarlagadda and John E. Hershey},
  doi       = {10.1007/978-1-4615-6313-6},
  year      = 1997,
  publisher = {Springer}
}

@book{Seberry2017,
  author  = {Jennifer Seberry},
  doi     = {10.1007/978-3-319-59032-5},
  journal = {Orthogonal Designs: Hadamard Matrices, Quadratic Forms and Algebras},
  title   = {Orthogonal designs: Hadamard matrices, quadratic forms and algebras},
  year    = {2017}
}

@article{Kharaghani2005,
  author  = {H. Kharaghani and B. Tayfeh-Rezaie},
  doi     = {10.1002/jcd.20043},
  issn    = {1063-8539},
  issue   = {6},
  journal = {Journal of Combinatorial Designs},
  pages   = {435-440},
  title   = {A Hadamard matrix of order 428},
  volume  = {13},
  year    = {2005}
}

@article{Spence1975,
  author  = {Edward Spence},
  doi     = {10.1016/0097-3165(75)90054-0},
  issn    = {00973165},
  issue   = {3},
  journal = {Journal of Combinatorial Theory, Series A},
  pages   = {287-300},
  title   = {Hadamard matrices from relative difference sets},
  volume  = {19},
  year    = {1975}
}

@article{Spence1975b,
  author  = {Edward Spence},
  doi     = {10.4153/cjm-1975-066-9},
  issn    = {0008-414X},
  issue   = {3},
  journal = {Canadian Journal of Mathematics},
  pages   = {555-560},
  title   = {Skew-Hadamard Matrices of the Goethals-Seidel Type},
  volume  = {27},
  year    = {1975}
}

@article{Djokovic1992a,
  author  = {Dragomir Ž. Djoković},
  doi     = {10.1017/S0004972700030185},
  issn    = {0004-9727},
  issue   = {2},
  journal = {Bulletin of the Australian Mathematical Society},
  pages   = {327-332},
  title   = {Construction of some new Hadamard matrices},
  volume  = {45},
  year    = {1992}
}

@article{Djokovic1992b,
  author  = {Dragomir Ž. Djoković},
  doi     = {10.1016/0097-3165(92)90029-T},
  issn    = {00973165},
  issue   = {2},
  journal = {Journal of Combinatorial Theory, Series A},
  pages   = {319-321},
  title   = {Skew Hadamard matrices of order 4 × 37 and 4 × 43},
  volume  = {61},
  year    = {1992}
}

@article{Djokovic1992c,
  author  = {Dragomir Ž. Djoković},
  journal = {Publikacije Elektrotehničkog fakulteta. Serija Matematika},
  pages   = {47-59},
  title   = {Ten New Orders of Hadamard Matrices of Skew Type},
  volume  = {3},
  year    = {1992}
}

@article{Djokovic1994a,
  author  = {Dragomir Ž. Djoković},
  journal = {Australasian Journal of Combinatorics},
  pages   = {259-264},
  title   = {Five New Orders for Hadamard Matrices of Skew Type},
  volume  = {10},
  year    = {1994}
}

@article{Szekeres1969,
  author  = {George Szekeres},
  doi     = {10.5169/seals-43227},
  journal = {L'Enseignement Mathématique},
  title   = {Tournaments and Hadamard Matrices},
  volume  = {15},
  year    = {1969}
}

@article{Szekeres1971,
  author  = {G. Szekeres},
  doi     = {10.4064/aa-18-1-349-353},
  issn    = {0065-1036},
  journal = {Acta Arithmetica},
  pages   = {349-353},
  title   = {Cyclotomy and complementary difference sets},
  volume  = {18},
  year    = {1971}
}

@article{Cooper1972,
  author  = {Joan Cooper and Jennifer Wallis},
  doi     = {10.1017/S0004972700045019},
  issn    = {0004-9727},
  issue   = {2},
  journal = {Bulletin of the Australian Mathematical Society},
  pages   = {269-277},
  title   = {A construction for Hadamard arrays},
  volume  = {7},
  year    = {1972}
}

@article{Sawade1985,
  author  = {Kazue Sawade},
  doi     = {10.1007/BF02582942},
  issn    = {0911-0119},
  issue   = {1},
  journal = {Graphs and Combinatorics},
  pages   = {185-187},
  title   = {A Hadamard matrix of order 268},
  volume  = {1},
  year    = {1985}
}

@article{Szekeres1987,
  author  = {George Szekeres},
  journal = {Combinatorics. Proceedings of the 7th Hungarian colloquium},
  pages   = {489-498},
  title   = {A note on skew type orthogonal ±1  matrices},
  year    = {1987}
}

@article{Koukouvinos2008,
  author  = {C. Koukouvinos and S. Stylianou},
  doi     = {10.1016/j.disc.2006.06.037},
  issn    = {0012365X},
  issue   = {13},
  journal = {Discrete Mathematics},
  pages   = {2723-2731},
  title   = {On skew-Hadamard matrices},
  volume  = {308},
  year    = {2008}
}

@article{Miyamoto1991,
  author  = {Masahiko Miyamoto},
  doi     = {10.1016/0097-3165(91)90008-5},
  issn    = {00973165},
  issue   = {1},
  journal = {Journal of Combinatorial Theory, Series A},
  pages   = {86-108},
  title   = {A construction of Hadamard matrices},
  volume  = {57},
  year    = {1991}
}

@article{Blatt1969,
  author  = {D. Blatt and G. Szekeres},
  doi     = {10.4153/CJM-1969-144-2},
  issn    = {0008-414X},
  journal = {Canadian Journal of Mathematics},
  pages   = {1319-1322},
  title   = {A Skew Hadamard Matrix of Order 52},
  volume  = {21},
  year    = {1969}
}

@article{Seberry1978,
  author  = {Jennifer Seberry},
  journal = {Ars Combinatoria},
  pages   = {255-276},
  title   = {On Skew Hadamard Matrices},
  volume  = {6},
  year    = {1978}
}

@article{Seberry1991,
  author     = {Seberry, Jennifer and Yamada, Mieko},
  title      = {Amicable {H}adamard matrices and amicable orthogonal designs},
  journal    = {Utilitas Math.},
  fjournal   = {Utilitas Mathematica. An International Journal of Discrete and
                Combinatorial Mathematics, and Statistical Design},
  volume     = {40},
  year       = {1991},
  pages      = {179--192},
  issn       = {0315-3681},
  mrclass    = {05B20},
  mrnumber   = {1143268},
  mrreviewer = {R.\ C.\ Mullin}
}

@article{Hadamard1893,
  author  = {Jacques Hadamard},
  journal = {Bulletin des Sciences Mathematiques},
  pages   = {240-246},
  title   = {Résolution d’une question relative aux determinants},
  volume  = {17},
  year    = {1893}
}

@article{Djokovic2016,
  author  = {Dragomir Ž. Djoković and Ilias S. Kotsireas},
  doi     = {10.1515/spma-2016-0029},
  issn    = {2300-7451},
  issue   = {1},
  journal = {Special Matrices},
  pages   = {317-325},
  title   = {A class of cyclic $(v; k_1, k_2, k_3; \lambda)$ difference families with $v \equiv 3$ (mod 4) a prime},
  volume  = {4},
  year    = {2016}
}

@online{cati2023,
  title         = {Implementing Hadamard Matrices in SageMath},
  author        = {Matteo Cati and Dmitrii V. Pasechnik},
  year          = {2023},
  eprint        = {2306.16812},
  archiveprefix = {arXiv},
  primaryclass  = {math.CO}
}

@online{Djokovic2010,
  title         = {Small orders of Hadamard matrices and base sequences},
  author        = {Dragomir Ž. Djoković},
  year          = {2010},
  eprint        = {1008.2043},
  archiveprefix = {arXiv},
  primaryclass  = {math.CO}
}

@online{Djokovic2023pers,
  title  = {Personal communucation},
  author = {Dragomir Ž. Djoković},
  year   = {2023},
  month  = {1}
}

@online{djokovic2024two,
  title         = {Two classes of Hadamard matrices of Goethals-Seidel type},
  author        = {Dragomir Ž. Djoković},
  year          = {2024},
  eprint        = {2404.14375},
  archiveprefix = {arXiv},
  primaryclass  = {math.CO}
}

@article{Paley1933,
  author  = {R. E. A. C. Paley},
  doi     = {10.1002/sapm1933121311},
  issn    = {00971421},
  issue   = {1-4},
  journal = {Journal of Mathematics and Physics},
  pages   = {311-320},
  title   = {On Orthogonal Matrices},
  volume  = {12},
  year    = {1933}
}

@article{Sylvester1867,
  author  = {J.J. Sylvester},
  doi     = {10.1080/14786446708639914},
  issn    = {1941-5982},
  issue   = {232},
  journal = {The London, Edinburgh, and Dublin Philosophical Magazine and Journal of Science},
  pages   = {461-475},
  title   = {Thoughts on inverse orthogonal matrices, simultaneous signsuccessions, and tessellated pavements in two or more colours, with applications to Newton's rule, ornamental tile-work, and the theory of numbers},
  volume  = {34},
  year    = {1867}
}

@article{Goethals1970,
  author  = {J. M. Goethals and J. J. Seidel},
  doi     = {10.1017/S144678870000673X},
  issn    = {0004-9735},
  issue   = {3},
  journal = {Journal of the Australian Mathematical Society},
  pages   = {343-344},
  title   = {A skew Hadamard matrix of order 36},
  volume  = {11},
  year    = {1970}
}

@article{Wallis1971,
  author  = {Jennifer Wallis},
  doi     = {10.1017/S0004972700047079},
  issn    = {0004-9727},
  issue   = {2},
  journal = {Bulletin of the Australian Mathematical Society},
  pages   = {203-204},
  title   = {A skew-Hadamard matrix of order 92},
  volume  = {5},
  year    = {1971}
}

@article{Cohen2017,
  author  = {Nathann Cohen and Dmitrii V. Pasechnik},
  doi     = {10.1007/s10623-016-0264-x},
  issn    = {0925-1022},
  issue   = {1-2},
  journal = {Designs, Codes and Cryptography},
  pages   = {223-235},
  title   = {Implementing Brouwer’s database of strongly regular graphs},
  volume  = {84},
  year    = {2017}
}

@online{Djokovic2023,
  title         = {Skew-Hadamard matrices of order 276},
  author        = {Dragomir Ž. Djoković},
  year          = {2023},
  eprint        = {2301.02751},
  archiveprefix = {arXiv},
  primaryclass  = {math.CO}
}

@manual{sagemath,
  key    = {SageMath},
  author = {{The Sage Developers}},
  title  = {{S}ageMath, the {S}age {M}athematics {S}oftware {S}ystem ({V}ersion 10.0)},
  url    = {https://www.sagemath.org},
  year   = {2023}
}

@book{MR3187587,
  author    = {London, Stephen},
  title     = {\href{http://gateway.proquest.com/openurl?url_ver=Z39.88-2004&rft_val_fmt=info:ofi/fmt:kev:mtx:dissertation&res_dat=xri:pqm&rft_dat=xri:pqdiss:3573367}{Constructing {N}ew {T}uryn {T}ype {S}equences, {T}-{S}equences
               and {H}ada\-mard {M}atrices}},
  note      = {Thesis (Ph.D.)--University of Illinois at Chicago},
  publisher = {ProQuest LLC, Ann Arbor, MI},
  year      = {2013},
  pages     = {64},
  isbn      = {978-1303-43542-3},
  mrclass   = {99-05},
  mrnumber  = {3187587}
}

@book{Seberry2020,
  title     = {Hadamard Matrices: Constructions using number theory and linear algebra},
  author    = {Seberry, Jennifer and Yamada, Mieko},
  year      = {2020},
  publisher = {John Wiley \& Sons},
  isbn      = {9781119520245},
  doi       = {10.1002/9781119520252}
}

@online{aebdb,
  author = {Andries E. Brouwer},
  title  = {Parameters of Strongly Regular Graphs},
  url    = {https://www.win.tue.nl/~aeb/graphs/srg/srgtab.html}
}

@article{ddj2022,
  author  = {Dragomir Ž. Djoković},
  title   = {Some new symmetric Hadamard matrices},
  doi     = {10.31799/1684-8853-2022-2-2-10},
  number  = {2},
  journal = {Information and Control Systems},
  year    = {2022},
  month   = {5},
  pages   = {2-10}
}

@misc{luming2021,
  title        = {Formal Verification of Constructions and
                  Theorems on Hadamard Matrices},
  author       = {Zhang, Lu-Ming},
  year         = 2021,
  howpublished = {MSc dissertation},
  address      = {Univesity of Oxford}
}

@article{Suksmono2022,
  author   = {Suksmono, Andriyan Bayu
              and Minato, Yuichiro},
  title    = {Quantum computing formulation of some classical Hadamard matrix searching methods and its implementation on a quantum computer},
  journal  = {Nature Scientific Reports},
  year     = {2022},
  volume   = {12},
  number   = {1},
  pages    = {197},
  issn     = {2045-2322},
  doi      = {10.1038/s41598-021-03586-0}
}

@online{magmatables,
  author  = {Computational Algebra Group},
  address = {University of Sydney},
  year    = 2024,
  title   = {Magma: Hadamard Matrices},
  url     = {http://magma.maths.usyd.edu.au/magma/info/topic_cmb_thy_had}
}

@article{MR1484478,
  author   = {Bosma, Wieb and Cannon, John and Playoust, Catherine},
  title    = {The {M}agma algebra system. {I}. {T}he user language},
  note     = {Computational algebra and number theory (London, 1993)},
  journal  = {J. Symbolic Comput.},
  fjournal = {Journal of Symbolic Computation},
  volume   = {24},
  year     = {1997},
  number   = {3-4},
  pages    = {235--265},
  issn     = {0747-7171},
  mrclass  = {68Q40},
  mrnumber = {MR1484478},
  doi      = {10.1006/jsco.1996.0125}
}

@online{sloantables,
  author  = {N. J. A. Sloane},
  address = {online},
  title   = {A Library of Hadamard Matrices},
  url     = {http://neilsloane.com/hadamard/}
}

@online{Rhadamard,
  author  = {Thomas Lumley et. al},
  address = {online},
  title   = {R package ``survey'': Hadamard matrices},
  url     = {https://search.r-project.org/CRAN/refmans/survey/html/hadamard.html}
}

@manual{Rsystem,
  title        = {R: A Language and Environment for Statistical
                  Computing},
  author       = {{R Core Team}},
  organization = {R Foundation for Statistical Computing},
  address      = {Vienna, Austria},
  year         = 2024,
  url          = {https://www.R-project.org}
}

@article{Djokovic2008,
  author   = {Đokovi\'{c}, Dragomir \v{Z}.},
  title    = {Hadamard matrices of order 764 exist},
  journal  = {Combinatorica},
  fjournal = {Combinatorica. An International Journal on Combinatorics and
              the Theory of Computing},
  volume   = {28},
  year     = {2008},
  number   = {4},
  pages    = {487--489},
  issn     = {0209-9683,1439-6912},
  mrclass  = {05B20 (05B30)},
  mrnumber = {2452846},
  doi      = {10.1007/s00493-008-2384-z}
}

@online{oeisA040081,
  title = {Riesel problem},
  url   = {https://oeis.org/A040081}
}

@online{Rieselwiki,
  title = {Riesel number},
  url   = {https://en.wikipedia.org/wiki/Riesel_number}
}

@online{hadgen,
  title  = {On-line utility to obtain all orders up to 1000, except 668, 716, 876 \& 892.},
  author = {Gupta, V.K. and Parsad Rajender and  Dhandapani, A.},
  url    = {https://drs.icar.gov.in/WebHadamard/WebHadamard.htm},
  year   = 2007
}

@article{Riesel1956,
  author  = {Riesel, Hans},
  title   = {Några stora primtal},
  journal = {Elementa},
  year    = 1956,
  pages   = {258--260},
  volume  = 39
}

@book{RGuy,
  author     = {Guy, Richard K.},
  title      = {Unsolved problems in number theory},
  series     = {Problem Books in Mathematics},
  edition    = {Third},
  publisher  = {Springer-Verlag, New York},
  year       = {2004},
  pages      = {xviii+437},
  isbn       = {0-387-20860-7},
  mrclass    = {11-01 (11-02)},
  mrnumber   = {2076335},
  mrreviewer = {Greg\ Martin},
  doi        = {10.1007/978-0-387-26677-0}
}

@online{MORiesel,
  title  = {Solving $2^{x+1} m -1=p^y$ for prime $p$ and natural $x,y$},
  author = {Max Alekseyev},
  eprint = {https://mathoverflow.net/q/471406},
  year   = 2024
}

@article{DELAUNEY20091002,
  title    = {On the asymptotic existence of Hadamard matrices},
  journal  = {Journal of Combinatorial Theory, Series A},
  volume   = {116},
  number   = {4},
  pages    = {1002-1008},
  year     = {2009},
  issn     = {0097-3165},
  doi      = {10.1016/j.jcta.2009.01.001},
  author   = {Warwick {de Launey}}
}

@article{Mat1978,
  author  = {R. A. Mathon},
  title   = {Symmetric conference matrices of order $pq^2 + 1$},
  journal = {Canad. J. Math.},
  volume  = 30,
  year    = 1978,
  pages   = {321-331},
  doi     = {10.4153/CJM-1978-029-1}
}

@article{tropp2011,
  author   = {Tropp, Joel A.},
  Title    = {Improved Analysis of the Subsampled Randomized Hadamard Transform},
  journal  = {Advances in Adaptive Data Analysis},
  volume   = {03},
  number   = {1-2},
  pages    = {115-126},
  year     = {2011},
  doi      = {10.1142/S1793536911000787}
}

@inproceedings{davenport2023,
  title     = {Proving an Execution of an Algorithm Correct?},
  author    = {Davenport, James H.},
  year      = {2023},
  month     = aug,
  day       = {28},
  doi       = {10.1007/978-3-031-42753-4_17},
  language  = {English},
  isbn      = {9783031427527},
  series    = {Lecture Notes in Computer Science},
  publisher = {Springer Science and Business Media Deutschland GmbH},
  pages     = {255--269},
  editor    = {Catherine Dubois and Manfred Kerber},
  booktitle = {Intelligent Computer Mathematics - 16th International Conference, CICM 2023, Proceedings},
  address   = {Germany}
}

@article{bauer98,
  author   = {Bauer, Andrej and Clarke, Edmund and Zhao, Xudong},
  year     = 1998,
  title    = {Analytica --- An Experiment in Combining Theorem Proving and Symbolic Computation},
  journal  = {Journal of Automated Reasoning},
  pages    = {295-325},
  volume   = 21,
  doi      = {10.1023/A:1006079212546}
}

@eprint{lewis2021,
  title         = {A bi-directional extensible interface between Lean and Mathematica},
  author        = {Robert Y. Lewis and Minchao Wu},
  year          = {2021},
  eprint        = {2101.07758},
  archiveprefix = {arXiv},
  primaryclass  = {cs.LO},
  doi           = {10.48550/arXiv.2101.07758}
}

@article{Djokovic2008a,
title = {Skew-Hadamard matrices of orders 188 and 388 exist},
journal={International Mathematical Forum},
volume=3,
number=22,
pages = {1063-1068},
  author   = {Đokovi\'{c}, Dragomir \v{Z}.},
  year     = {2008},
  doi={10.48550/arXiv.0704.0640},
}

@article{Djokovic2008b,
author = {Đoković, Dragomir Ž.},
title = {Skew-Hadamard matrices of orders 436, 580, and 988 exist},
journal = {Journal of Combinatorial Designs},
volume = {16},
number = {6},
pages = {493-498},
doi = {10.1002/jcd.20180},
year = {2008},
}

@article{DGK2014,
author = {Đoković, Dragomir Ž. and Golubitsky, Oleg and Kotsireas, Ilias S.},
title = {Some New Orders of Hadamard and Skew-Hadamard Matrices},
journal = {Journal of Combinatorial Designs},
volume = {22},
number = {6},
pages = {270-277},
doi = {doi.org/10.1002/jcd.21358},
year = {2014},
}

\section*{Appendix: Tables of constructions and known orders}
\label{appendix}
Table \ref{hadamard-constructions-list-table} details for every odd value of $n < 300$ which method can be used to construct the Hadamard matrix of order $4n$ (whenever $n$ is even, the Hadamard matrix can be constructed using the doubling construction). Similarly, Table \ref{skew-hadamard-constructions-list-table} lists the constructions of skew Hadamard matrices. Note that some entries of the two tables are empty: this indicates that no construction for the corresponding order is known.
Table \ref{construction-list-table} contains an explanation of the abbreviations used in the two lists. 

\begin{table}[hbt]
	\caption{Algorithms used in Tables \ref{hadamard-constructions-list-table}-\ref{skew-hadamard-constructions-list-table}.
\label{construction-list-table}}
\begin{centering}
\begin{tabular}{>{\raggedright\arraybackslash}p{0.15\linewidth} | >{\raggedright\arraybackslash}p{0.72\linewidth} |> {\raggedright\arraybackslash}p{0.10\linewidth} }
 PaleyI &  Paley's first construction & Sect. \ref{paleyconst}\\
 PaleyII & Paley's second construction & Sect. \ref{paleyconst}\\
 Will &  Williamson construction & Sect. \ref{sect:willi} \\
 GS & Goethals-Seidel array  & Sect. \ref{gsarray} \\
 SDS & Construction from (possibly skew) SDS  & Sect. \ref{sect:SDS} \\
 CW(t) & Construction of Hadamard matrix of order $4n$ from T-sequ\-ences of length $t$ (and Williamson matrices of order $n/t$)  & Sect. \ref{sect:Tmat} \\
 Good & construction from good matrices  & Sect. \ref{sect:good} \\
 Miy & Miyamoto Construction  & Thm. \ref{miyamoto} \\
 CDS & Construction from complementary difference sets  & Sect. \ref{sect:SDS} \\
 Spence(q) & Spence construction of skew Hadamard matrix  of order $4(1+q+q^2)$  & Thm. \ref{thm:spence75b} \\
 AOD(m, n) & construction of a skew Hadamard matrix of order $mn(n-1)$ from amicable orthogonal designs  & Thm. \ref{thm:amicable} \\
 \end{tabular}
 \end{centering}
\end{table}


\begin{table}[hbt]
    \caption{Hadamard matrices of order $4n$, $n$ odd, up to 1200
    \label{hadamard-constructions-list-table}}
    \centering
    \begin{tabular}{|c c || c c || c c || c c || c c|}
    \hline 
        1 & PaleyI &
        3 & PaleyI &
        5 & PaleyI &
        7 & PaleyI &
        9 & PaleyII \\
        11 & PaleyI &
        13 & PaleyII &
        15 & PaleyI &
        17 & PaleyI &
        19 & PaleyII \\
        21 & PaleyI &
        23 & Will &
        25 & PaleyII &
        27 & PaleyI &
        29 & Will \\
        31 & PaleyII &
        33 & PaleyI &
        35 & PaleyI &
        37 & PaleyII &
        39 & Will \\
        41 & PaleyI &
        43 & Will &
        45 & PaleyI &
        47 & CW(47) &
        49 & PaleyII \\
        51 & PaleyII &
        53 & PaleyI &
        55 & PaleyII &
        57 & PaleyI &
        59 & CW(59) \\
        61 & PaleyI &
        63 & PaleyI &
        65 & CW(5) &
        67 & CW(67) &
        69 & PaleyII \\
        71 & PaleyI &
        73 & Miy &
        75 & PaleyII &
        77 & PaleyI &
        79 & PaleyII \\
        81 & CW(3) &
        83 & PaleyI &
        85 & PaleyII &
        87 & PaleyI &
        89 & CW(89) \\
        91 & PaleyII &
        93 & CW(3) &
        95 & PaleyI &
        97 & PaleyII &
        99 & PaleyII \\
        101 & Miy &
        103 & SDS &
        105 & PaleyI &
        107 & CW(107) &
        109 & Miy \\
        111 & PaleyI &
        113 & Miy &
        115 & PaleyII &
        117 & PaleyI &
        119 & CW(7) \\
        121 & PaleyII &
        123 & PaleyI &
        125 & PaleyI &
        127 & SDS &
        129 & PaleyII \\
        131 & PaleyI &
        133 & CW(7) &
        135 & PaleyII &
        137 & PaleyI &
        139 & PaleyII \\
        141 & PaleyI &
        143 & PaleyI &
        145 & PaleyII &
        147 & PaleyI &
        149 & Miy \\
        151 & SDS &
        153 & CW(3) &
        155 & PaleyI &
        157 & PaleyII &
        159 & PaleyII \\
        161 & PaleyI &
        163 & SDS &
        165 & PaleyI &
        167 &  &
        169 & PaleyII \\
        171 & PaleyI &
        173 & PaleyI &
        175 & PaleyII &
        177 & PaleyII &
        179 &  \\
        181 & PaleyII &
        183 & CW(3) &
        185 & PaleyI &
        187 & PaleyII &
        189 & CW(3) \\
        191 & SDS &
        193 & Miy &
        195 & PaleyII &
        197 & PaleyI &
        199 & PaleyII \\
        201 & PaleyII &
        203 & PaleyI &
        205 & PaleyII &
        207 & PaleyI &
        209 & CW(11) \\
        211 & PaleyII &
        213 & CW(71) &
        215 & PaleyI &
        217 & PaleyII &
        219 & SDS \\
        221 & PaleyI &
        223 &  &
        225 & PaleyII &
        227 & PaleyI &
        229 & PaleyII \\
        231 & PaleyII &
        233 & Miy &
        235 & CW(47) &
        237 & PaleyI &
        239 & SDS \\
        241 & Miy &
        243 & PaleyI &
        245 & CW(5) &
        247 & CW(13) &
        249 & CW(83)  \\
        251 & SDS &
        253 & CW(11) &
        255 & PaleyI &
        257 & Miy &
        259 & CW(7)  \\
        261 & PaleyII  &
        263 & PaleyI &
        265 & PaleyII &
        267 & CW(89) &
        269 & Miy \\
        271 & PaleyII &
        273 & PaleyI &
        275 & CW(5) &
        277 & Miy &
        279 & PaleyII  \\
        281 & PaleyI &
        283 &  &
        285 & PaleyII &
        287 & CW(7) &
        289 & PaleyII \\
        291 & PaleyI &
        293 & PaleyI &
        295 & CW(59) &
        297 & PaleyI &
        299 & CW(13) \\
        301 & PaleyII \\
    \hline 
    \end{tabular}
\end{table}

\begin{table}[hbt]
    \caption{Skew Hadamard matrices of order $4n$, $n$ odd, up to 1200
    \label{skew-hadamard-constructions-list-table}}
    \centering
    \begin{tabular}{|c c || c c || c c || c c || c c|}
    \hline 
        1 & Good &
        3 & Good &
        5 & Good &
        7 & Good &
        9 & Good \\
        11 & Good &
        13 & Good &
        15 & Good &
        17 & Good &
        19 & Good \\
        21 & Good &
        23 & Good &
        25 & Good &
        27 & Good &
        29 & Good \\
        31 & Good &
        33 & PaleyI &
        35 & PaleyI &
        37 & SDS &
        39 & SDS \\
        41 & PaleyI &
        43 & SDS &
        45 & PaleyI &
        47 & GS &
        49 & SDS \\
        51 & CDS &
        53 & PaleyI &
        55 & CDS &
        57 & PaleyI &
        59 & GS \\
        61 & PaleyI &
        63 & PaleyI &
        65 & SDS &
        67 & SDS &
        69 & GS \\
        71 & PaleyI &
        73 & SDS &
        75 & CDS &
        77 & PaleyI &
        79 & CDS \\
        81 & SDS &
        83 & PaleyI &
        85 & CDS &
        87 & PaleyI &
        89 &   \\
        91 & CDS &
        93 & SDS &
        95 & PaleyI &
        97 & SDS &
        99 & CDS \\
        101 &   &
        103 & SDS &
        105 & PaleyI &
        107 &   &
        109 & SDS \\
        111 & PaleyI &
        113 & SDS &
        115 & CDS &
        117 & PaleyI &
        119 &   \\
        121 & SDS &
        123 & PaleyI &
        125 & PaleyI &
        127 & SDS &
        129 & SDS \\
        131 & PaleyI &
        133 & SDS &
        135 & CDS &
        137 & PaleyI &
        139 & CDS \\
        141 & PaleyI &
        143 & PaleyI &
        145 & SDS &
        147 & PaleyI &
        149 &   \\
        151 & SDS &
        153 &   &
        155 & PaleyI &
        157 & SDS &
        159 & CDS \\
        161 & PaleyI &
        163 & SDS &
        165 & PaleyI &
        167 &   &
        169 & SDS \\
        171 & PaleyI &
        173 & PaleyI &
        175 & CDS &
        177 &   &
        179 &   \\
        181 & SDS &
        183 & Spence(13) &
        185 & PaleyI &
        187 & CDS &
        189 & AOD(1, 28) \\
        191 &   &
        193 &   &
        195 & CDS &
        197 & PaleyI &
        199 & CDS \\
        201 &   &
        203 & PaleyI &
        205 &   &
        207 & PaleyI &
        209 &   \\
        211 & CDS &
        213 & SDS &
        215 & PaleyI &
        217 & SDS &
        219 & SDS \\
        221 & PaleyI &
        223 &   &
        225 &   &
        227 & PaleyI &
        229 &   \\
        231 & CDS &
        233 &   &
        235 &   &
        237 & PaleyI &
        239 & SDS \\
        241 & SDS &
        243 & PaleyI &
        245 &   &
        247 & SDS &
        249 &   \\
        251 &  &
        253 & &
        255 & PaleyI &
        257 & &
        259 & \\
        261 & &
        263 &  PaleyI &
        265 & &
        267 & SDS &
        269 & \\
        271 & CDS &
        273 & PaleyI &
        275 & &
        277 & &
        279 & CDS \\
        281 & PaleyI &
        283  & &
        285 & &
        287 & &
        289 & \\
        291 &  PaleyI &
        293 &  PaleyI &
        295 & &
        297 & PaleyI &
        299 & \\
    \hline 
    \end{tabular}
\end{table}

\afterpage{
\begin{table}[H]
    \caption{Pairs $n(m)$ where $m$ is minimal, s.t. a Hadamard matrix of order $2^mn$ is known,  $n \le 2999$.
    The case $m=2$ is not shown.\protect\footnotemark 
More details in Sect.~\ref{knownhm}.
    \label{hadamard-known-list-table}
	}
    \centering
    \begin{tabular}{|l l l l l l l l|}
    \hline 
167(3) & 179(3) & 223(3) & 283(3) & 311(3) & 347(3) & 359(4) & 419(3) \\
443(3) & 479(4) & 487(3) & 491(5) & 515(3) & 523(3) & 537(3) & 571(3) \\
573(3) & 599(6) & 643(3) & 647(3) & 659(4) & 669(3) & 719(4) & 721(3) \\
739(6) & 751(3) & 789(3) & 823(3) & 839(4) & 859(3) & 863(3) & 883(3) \\
907(3)  & 917(3)  & 919(3)  & 933(4)  & 947(3)  & 955(3)  & 971(6)  & 991(3) \\
1019(3) & 1031(6) & 1039(3) & 1051(3) & 1063(3) & 1087(3) & 1103(3) & 1115(3) \\
1123(3) & 1133(3) & 1169(4) & 1187(3) & 1223(6) & 1255(3) & 1257(3) & 1259(4) \\
1283(3) & 1291(3) & 1303(3) & 1315(3) & 1319(4) & 1327(3) & 1359(3) & 1367(3) \\
1423(3) & 1427(3) & 1437(6) & 1439(8) & 1441(3) & 1447(7) & 1451(6) & 1471(3) \\
1473(3) & 1483(3) & 1487(3) & 1499(7) & 1509(3) & 1527(3) & 1543(3) & 1559(4) \\
1567(6) & 1571(7) & 1579(5) & 1583(3) & 1589(3) & 1619(3) & 1661(3) & 1663(3) \\
1667(3) & 1689(3) & 1699(3) & 1703(3) & 1713(3) & 1719(3) & 1747(3) & 1751(3) \\
1783(6) & 1787(3) & 1793(4) & 1795(5) & 1823(3) & 1831(3) & 1841(3) & 1847(3) \\
1871(3) & 1879(3) & 1883(3) & 1893(3) & 1907(3) & 1915(3) & 1929(3) & 1949(4) \\
1957(3) & 1963(3) & 1969(4) & 1979(4) & 1981(4) & 1987(7) & 2027(3) & 2039(7) \\
2063(7) & 2083(3) & 2087(4) & 2095(3) & 2099(3) & 2119(4) & 2143(3) & 2155(3) \\
2171(4) & 2203(3) & 2207(4) & 2215(4) & 2227(3) & 2251(5) & 2287(7) & 2293(6) \\
2315(4) & 2327(4) & 2335(3) & 2339(4) & 2347(3) & 2369(3) & 2371(6) & 2383(3) \\
2399(8) & 2423(4) & 2429(4) & 2459(3) & 2489(3) & 2503(3) & 2513(4) & 2515(3) \\
2531(3) & 2543(6) & 2545(3) & 2571(3) & 2579(3) & 2589(4) & 2591(3) & 2629(3) \\
2647(3) & 2659(3) & 2661(3) & 2671(7) & 2677(6) & 2683(7) & 2687(9) & 2699(7) \\
2711(3) & 2731(3) & 2733(3) & 2767(3) & 2803(3) & 2815(3) & 2819(3) & 2823(3) \\
2841(3) & 2843(3) & 2855(4) & 2865(3) & 2879(7) & 2887(5) & 2893(3) & 2899(4) \\
2903(4) & 2913(6) & 2927(3) & 2939(6) & 2951(3) & 2963(3) & 2971(3) & 2973(4) \\
2987(3) & 2995(6) & 2999(9) &         &         &         &         & \\
    \hline 
 \end{tabular}
\end{table}

    \footnotetext{The following values are in \cite{Seberry2020}, Table A.17, with $m>2$, but not here:
    191 - cf. Sect.~\ref{gsarray} \cite{Djokovic2008}; 
    787 - cf. Theorem~\ref{paleyII-theo} with a symmetric conference matrix from \cite{Mat1978};
    749, 781, 1177, 1349, 1491, 1633, 2033, 2059, 2461, 2627, 2769 - cf. Theorem~\ref{CooperWallisConstruction-section} \cite{Cooper1972}; 
    853, 1093, 1669, 1789, 1913, 1933, 2053, 2269, 2341 - cf. Theorem~\ref{miyamoto} \cite{Miyamoto1991}; 1723 - cf. Sect.~\ref{skewmats} \cite{Spence1975b}.}
}

\afterpage{
\begin{table}
    \centering
    \begin{tabular}{|l l l l l l l l l|}
    \hline 
89(4) & 101(10) & \rb{107(10)} & 119(4) & 149(4) & 153(3) & 167(4) & \pb{177(4)} & 179(8) \\
\rb{191(226)} & 193(3) & 201(3) & 205(3) & 209(4) & 223(3) & 225(4) & 229(3) & 233(4) \\
235(3) & 245(4) & 249(4) & 251(6) & 253(4) & 257(4) & 259(4) & 261(3) & 265(4) \\
269(8) & 275(4) & 277(5) & \rb{283(11)} & 285(3) & 287(4) & 289(3) & \pb{295(4)} & 299(4) \\
301(3) & 303(3) & 305(4) & 309(3) & 311(26) & 317(6) & 319(3) & \pb{325(4)} & \pb{329(4)} \\
\pb{335(4)} & 337(18) & 341(4) & \pb{343(4)} & 345(4) & \rb{347(18)} & 349(3) & 353(4) & 359(4) \\
361(3) & 369(4) & 373(7) & \pb{377(5)} & 385(3) & 389(15) & 391(4) &  \rb{401(10)}  & 403(5) \\
409(3) & 413(4) & 419(4) & 423(4) & 429(3) & 433(3) & 435(4) & 441(3) & 443(6) \\
 445(3) & \rb{449(112)} & 455(4) & \rb{457(9)} & 459(3) & \rb{461(10)} & 463(7) & 465(3) & 469(3) \\
\pb{473(4)} & 475(4) & \rb{479(12)} & 481(3) & 485(4) & 487(5) & 489(3) & \rb{491(46)}  & 493(3) \\
\rb{505(11)} & \rb{509(8)} & 513(4) & \pb{515(4)} & \pb{517(4)} & 519(4) & 523(7) & 527(4) & 529(3) \\
537(5) & 539(4) & 541(3) & \pb{543(4)} & 549(3) & 553(3) & \rb{557(8)} & 565(3) & 569(4) \\
571(3) & 573(3) & 575(4) & \rb{577(29)} & 579(5) & 581(4) & 583(3) & 589(5) & 595(3) \\
597(4) & \rb{599(8)} & 601(5) & 605(4) & 609(3) & \pb{611(5)} & 613(3) & 621(3) & 623(4) \\
625(3) & 627(4) & \pb{629(4)} & 637(4) & 641(6) & \rb{643(47)} & 647(10) & \pb{649(4)} & 653(10) \\
655(4) & \pb{657(4)} & \rb{659(\scriptsize 800516)} & \rb{661(9)} & 667(6) & 669(3) & 673(7) & 679(3) & 681(4) \\
689(4) & 693(4) & 695(4) & 697(4) & 699(3) & 703(3) & 707(4) & \rb{709(69)} & 717(4) \\
719(4) & \pb{721(4)} & 723(3) & \pb{725(4)} & 729(4) & \pb{731(4)} & \rb{733(19)} & \pb{737(4)} & \rb{739(107)} \\
745(6) & 749(12) & 751(3) & 759(4) & 761(22) & \pb{763(4)} & 765(4) & 769(3) & \rb{773(8)} \\
777(4) & 779(4) & 781(3) & \pb{785(4)} & 787(5) & 789(3) & 793(3) & 795(3) & \pb{803(4)} \\
805(4) & 809(12) & 817(5) & 821(6) & 823(3) & \rb{829(53)} & \rb{839(8)} &\pb{ 841(6)} & \pb{845(4)} \\
849(3) & \pb{851(6)} & 853(3) & 857(4) & 859(3) & 861(4) & 863(4) & 865(4) & 869(4) \\
\rb{877(13)} & 881(6) & \rb{883(7)} & \pb{889(4)} & 891(3) & 897(5) & 899(6) & 901(3) & 903(4) \\
905(4) & 907(5) & 909(4) & 913(4) & 917(4) & 919(3) & 925(3) & 927(4) & \rb{929(12)} \\
\rb{933(10)} & 937(5) & 941(6) & 943(4) & 947(6) & 949(3) & \rb{953(14)} & 955(3) & 957(4) \\
959(4) & 961(3) & 965(4) & 969(4) & 971(6) & 973(4) & 979(5) & 985(3) & 989(4) \\
991(3) & \nottt{995(4)} & \rb{997(9)} & & & & &\\
    \hline 
 \end{tabular}
    \caption{Pairs $n(m)$ where $m$ is minimal, s.t. a skew Hadamard matrix of order $2^mn$ is known,
    $n \le 999$. The case $m=2$ is not shown.\protect\footnotemark
    Our updates are shown in font/colour: {\rb{san~serif/red}} for Paley construction, 
    {\pb{slanted/blue}} for Williamson's construction,  Theorem~\ref{lemma1415hall}.
More details in Sect.~\ref{knownhm}.
    \label{skew-hadamard-known-list-table} }
\end{table}

\footnotetext{
    The following values are in \cite{Seberry2020}, Table 9.2, with $m>2$, but not here:
    39, 47, 49, 65, 69, 81, 93, 97, 103, 109, 121, 129, 133, 145, 151, 169, 217, 219, 239, 247, 267, 331, 885 -
   cf.  Sect.~\ref{gsarray}, 331 - cf. Theorem~\ref{sect:spence}, 397 - cf. Sect.~\ref{gsarray} and \cite{djokovic2024two},  885 - cf. Sect.~\ref{paleyconst}.
}
}
\end{document}